\documentclass[12]{article}

\usepackage{fleqn}
\usepackage{graphicx,epsf,psfrag}
\usepackage{amsthm,amsfonts,amssymb,amsmath}

\begin{document}
\Large
\begin{center}
{\bf Veldkamp Spaces:\\ From (Dynkin) Diagrams to (Pauli) Groups}
\end{center}
\vspace*{-.1cm}
\large
\begin{center}
Metod Saniga,$^{1,2}$ Fr\'ed\'eric Holweck$^{1}$ and Petr Pracna$^{3}$ 
\end{center}
\vspace*{-.4cm}
\normalsize
\begin{center}

$^{1}$IRTES-M3M/UTBM, Universit\'e de Bourgogne Franche-Comt\'e\\ 
Campus S\'evenans, F-90010 Belfort Cedex, France\\ (metod.saniga@utbm.fr, frederic.holweck@utbm.fr)

\vspace*{.2cm}

$^{2}$Astronomical Institute, Slovak Academy of Sciences\\
SK-05960 Tatransk\' a Lomnica, Slovak Republic\\
(msaniga@astro.sk) 

\vspace*{.2cm}

$^{3}$J. Heyrovsk\' y Institute of Physical Chemistry, v.\,v.\,i.\\  Academy of Sciences of the Czech Republic \\ Dolej\v skova 3, CZ-18223 Prague, Czech Republic\\
(ppracna@seznam.cz)

\end{center}

\vspace*{-.3cm} \noindent \hrulefill

\vspace*{-.0cm} \noindent {\bf Abstract}

\noindent Regarding a Dynkin diagram as a specific point-line incidence structure (where each line has just two points), one can associate with it a Veldkamp space. Focusing on extended Dynkin diagrams of type $\widetilde{D}_n$, $4 \leq n \leq 8$, it is shown that the corresponding Veldkamp space always contains a distinguished copy of the projective space PG$(3,2)$. Proper labelling of the vertices of the diagram (for $4 \leq n \leq 7$) by particular elements of the two-qubit Pauli group establishes a bijection between the 15 elements of the group and the 15 points of the PG$(3,2)$. The bijection is such that the product of three elements lying on the same line is the identity and one also readily singles out that particular copy of the symplectic polar space $W(3,2)$ of the PG$(3,2)$ whose lines correspond to triples of mutually commuting elements of the group; in the latter case, in addition, we arrive at a unique copy of the Mermin-Peres magic square. In the case of $n=8$, a more natural labeling is that in terms of elements of the three-qubit Pauli group, furnishing a bijection between the 63 elements of the group and the 63 points of PG$(5,2)$, the latter being the maximum projective subspace of the corresponding Veldkamp space; here, the points of the distinguished PG$(3,2)$ are in a bijection with the elements of a two-qubit subgroup of the three-qubit Pauli group, yielding a three-qubit version of the Mermin-Peres square. Moreover, save for $n=4$, each Veldkamp space is also endowed with some `exceptional' point(s). Interestingly, two such points in the $n=8$ case define a unique Fano plane whose inherited three-qubit labels feature solely the Pauli matrix $Y$.

\vspace*{.3cm}


\vspace*{-.2cm} \noindent \hrulefill


\section{Introduction} 

Although the fact that the properties of certain finite groups and the structure of
certain finite geometries/point-line incidence structures are tied very closely to each other has been known in the mathematics community for a relatively long time, it was only some ten years ago when this fact was also recognized by physicists. There exists, in
particular, a large family of groups relevant for quantum information theory where  (non)commutativity of two
distinct elements can be expressed in the language of finite symplectic polar spaces and/or finite projective ring
lines. More recently, this link has also been employed to get deeper insights into the so-called black-hole-qubit
correspondence --- a still puzzling relation between the entropy of certain stringy black holes and the
entanglement properties of some small-level quantum systems. A concept that played a particular, yet rather implicit, role in the
latter context turned out to be that of the {\it Veldkamp space} of a point-line incidence structure \cite{buec}.

The relevance of the concept of Veldkamp space for (quantum) physics was first recognized in the context of the geometry of the generalized two-qubit Pauli group \cite{twoq}. It was demonstrated in detail that the corresponding Veldkamp space features three distinct
types of points, each having a distinguished physical meaning, namely: a maximum set of pairwise
non-commuting group elements, a set of six elements commuting with the given one, and a set of nine
elements forming the so-called Mermin-Peres `magic' square. An intriguing novelty, stemming from the structure of Veldkamp
lines, was the recognition of (uni- and tri-centric) triads and specific pentads of elements in addition to
the above-mentioned `classical' subsets. Based on these findings, Vrana and L\'evay \cite{vrle} were able to ascertain the structure of the Veldkamp space of the $N$-qubit Pauli group for an arbitrary $N$,
singling out appropriate point-line incidence structures within the associated symplectic polar spaces of rank $N$
and order two. 

Quantum contextuality is another important aspect of quantum information theory where the notion of Veldkamp
space entered the game in an essential way. Employing the structure of the split Cayley hexagon of order two,
the smallest non-trivial generalized hexagon and a distinguished subgeometry of the symplectic polar space
$W(5,2)$ of the three-qubit Pauli group, two of the authors and their colleagues \cite{sppl} got an intriguing finite-geometric
insight into the nature of a couple of `magic' three-qubit configurations proposed by Waegell and Aravind \cite{waar}, as either of the two was found to be uniquely extendible into a
geometric hyperplane (i.\,e., a Veldkamp point) of the hexagon, this being of a very special type. Moreover, an automorphism of order
seven of the hexagon gave birth, for either of the two, to six more replicas, each having the same “magical”
nature as the parent one. As one of the most symmetric three-qubit `magic configurations is the so-called
Mermin pentagram, these observations prompted Planat and two of the authors \cite{psh}  to have a closer look at automorphisms of the split Cayley hexagon of order two and find out that $W(5,2)$ contains
altogether 12\,096 copies of such a pentagram, this number being -- surprisingly -- equal to the order of the
automorphism group of the hexagon. In addition, the authors succeeded in singling out those types of points of
the Veldkamp space of the hexagon that contain pentagrams and observed that the number of Veldkamp points
of a given type times the number of pentagrams the particular Veldkamp point contains is almost always a
multiple of 12\,096.

In light of the above-given observations, it is not surprising that the split Cayley hexagon and its Veldkamp space also
occur in the context of the already mentioned black-hole-qubit correspondence. In particular, the
$PSL_2(7)$ subgroup of the automorphism group of the hexagon and its associated Coxeter subgeometry were
found to be intricately related to the $E_7$-symmetric black-hole entropy formula in string theory, given a prominent
role played by Veldkamp points that answer to Klein quadrics in the ambient projective space PG$(5,2)$ 
\cite{lsv}. On the other hand, the $E_6$-symmetric entropy formula
describing black holes and black strings in five dimensions is underlaid by the geometry of generalized quadrangle GQ$(2,4)$ and its Veldkamp space \cite{lsvp}; here, the two
pronounced truncations of the entropy formula with 15 and 11 charges correspond exactly to two distinct types of
Veldkamp points of GQ$(2,4)$ \cite{sglpv}.

Apart from these interesting physical applications, the notion of Veldkamp space 
has also been successfully used in a couple of purely mathematical contexts. On the one hand \cite{shhpp}, it helped us to ascertain finer traits of the nested structure of Segre varieties that are $N$-fold direct product of projective lines of size three, $S_{(N)} \equiv$ PG$(1,2) \times$ PG$(1,2) \times \cdots \times$ PG$(1,2)$, for the cases $2 \leq N \leq 4$. In particular, given the fact that $S_{(N)} =$ PG$(1,2) \times S_{(N-1)}$, a powerful diagrammatical recipe was found that shows how to fully recover the properties of Veldkamp points (i.\,e. geometric hyperplanes) of $S_{(N)}$ once we know the types (and cardinalities thereof) of Veldkamp lines of $S_{(N-1)}$ \cite{grsa}. On the other hand 
\cite{sahopr}, it led to a better understanding of an intriguing finite-geometrical underpinning of the multiplication tables of real Cayley-Dickson algebras $A_N$, for $3 \leq N \leq 6$, that admits generalization to any higher-dimensional $A_N$. The multiplication properties of imaginary units of the algebra $A_N$ are encoded in the structure of the projective space PG$(N-1,2)$ that is endowed with a refined structure stemming from particular properties of triples of imaginary units forming its lines. 
The concept of Veldkamp space was here invoked to account for this refinement, with the relevant point-line incidence structure being a binomial $\left({N+1 \choose 2}_{N-1}, {N+1 \choose 3}_{3}\right)$-configuration, or, equivalently, a combinatorial Grassmannian of type $G_2(N+1)$.  

The present paper offers another interesting set of examples illustrating the usefulness of the notion of Veldkamp space in a broader context.  
Namely, we shall start with a sequence of extended Dynkin diagrams of type $\widetilde{D}_n$, $4 \leq n \leq 8$, then consider each diagram as a point-line incidence structure, next look at this structure through projective subspaces of its Veldkamp space in order, after appropriate labeling of the vertices of the diagrams by elements of two- (or three-)qubit Pauli groups,  to recapture the well-known geometrical representations of these groups.

\section{Basic Concepts and Notation}

In this section we shall give a brief summary of basic concepts, symbols and notation employed in the sequel.

We start with a {\it point-line incidence structure} $\mathcal{C} = (\mathcal{P},\mathcal{L},I)$ where $\mathcal{P}$ and $\mathcal{L}$ are, respectively, sets of points and lines and where incidence $I \subseteq \mathcal{P} \times \mathcal{L}$ is a binary relation indicating which point-line pairs are incident (see, e.\,g., \cite{shult}).  The dual of a point-line incidence structure is the structure with points
and lines exchanged, and with the reversed incidence relation. In what follows we shall encounter only specific point-line incidence structures where every line has just two points and any two distinct points are joined by at most one line. 
The {\it order} of a point of $\mathcal{C}$ is the number of lines passing through it.
A {\it geometric hyperplane} of $\mathcal{C} = (\mathcal{P},\mathcal{L},I)$ is a proper subset of $\mathcal{P}$ such that a line from $\mathcal{C}$ either lies fully in the subset, or shares with it only one point.
 If $\mathcal{C}$ possesses geometric hyperplanes, then one can define the {\it Veldkamp space} of $\mathcal{C}$ as follows \cite{buec}: (i) a point of the Veldkamp space (also called a Veldkamp point of $\mathcal{C}$) is a geometric hyperplane $H$ of  $\mathcal{C}$
and (ii) a line of the Veldkamp space (also called a Veldkamp line of $\mathcal{C}$) is the collection $H'H''$ of all geometric hyperplanes $H$ of $\mathcal{C}$  such that $H' \cap H'' = H' \cap H = H'' \cap H$ or $H = H', H''$, where $H'$ and $H''$ are distinct geometric hyperplanes. 
As Veldkamp lines of a $\mathcal{C}$ with all lines of size two are found to be of different cardinalities, we shall focus only on those subgeometries
of the corresponding Veldkamp space whose lines are all of size {\it three}, being of the form $\{H', H'', \overline{H' \Delta H''}\}$; here, the symbol $\Delta$ stands for the symmetric difference of the two geometric hyperplanes and an overbar denotes the complement of the object indicated.

Further, let $V(d+1,q)$, $d \geq 1$, denote a rank-$(d+1)$ vector space over the Galois field GF$(q)$, $q$ being a power of a prime. Associated with this vector space is a $d$-dimensional {\it projective space} over GF$(q)$, PG$(d,q)$, whose points, lines, planes,$\ldots$, hyperplanes are rank-one, rank-two, rank-three,$\ldots$, rank-$d$ subspaces of $V(d+1,q)$; for $q=2$, this projective space features
$2^{d+1} - 1$ points and $(2^{d+1} - 1)(2^{d} - 1)/3$ lines (see, e.\,g., \cite{hita}). 
Given a PG$(2N-1,q)$ that is endowed with a non-degenerate symplectic form, the {\it symplectic polar space} $W(2N-1,q)$ in PG$(2N-1,q)$ is the space of all totally isotropic subspaces with respect to the non-degenerate symplectic form \cite{cam}, with its maximal totally isotropic subspaces, also
called {\it generators}, having dimension $N - 1$.  For $q=2$ this
polar space contains $|$PG$(2N-1, 2)| = 2^{2N} - 1 = 4^{N} - 1$
points and $(2+1)(2^2+1)\cdots(2^N+1)$ generators. 

Next, we need to mention generalized (complex) $N$-qubit Pauli groups (see, e.\,g., \cite{nc}), ${\cal P}_N$,  generated by $N$-fold tensor products of the matrices
\begin{eqnarray*}
I = \left(
\begin{array}{cc}
1 & 0 \\
0 & 1 \\
\end{array}
\right),~
X = \left(
\begin{array}{cc}
0 & 1 \\
1 & 0 \\
\end{array}
\right),~
Y = \left(
\begin{array}{cr}
0 & -i \\
i & 0 \\
\end{array}
\right)
~{\rm and}~
Z = \left(
\begin{array}{cr}
1 & 0 \\
0 & -1 \\
\end{array}
\right).
\end{eqnarray*}
Explicitly,
\begin{equation*}
{\cal P}_N = \{ i^{\alpha} A_1 \otimes A_2 \otimes\cdots\otimes A_N:~ A_k \in \{I, X, Y, Z \},~ k = 1, 2,\cdots,N,~\alpha \in \{0,1,2,3\} \}.
\end{equation*}
Here, we are only interested in their factor versions $\overline{{\cal P}}_N \equiv {\cal P}_N/{\cal Z}({\cal P}_N)$, where ${\cal Z}({\cal P}_N) = \{\pm I_{(1)} \otimes I_{(2)} \otimes \cdots \otimes I_{(N)}, \pm i I_{(1)} \otimes I_{(2)} \otimes \cdots \otimes I_{(N)}\}$. For a particular value of $N$, the $4^N - 1$ elements of $\overline{{\cal P}}_N \backslash \{I_{(1)} \otimes I_{(2)} \otimes \cdots \otimes I_{(N)}\}$  can bijectively be identified with the same number of points of
$W(2N-1, 2)$ in such a way that two commuting elements of the group will lie on the same totally isotropic line and a maximum set of mutually commuting elements corresponds to a generator of this symplectic space (see, e.\,g., \cite{sp}--\cite{th}). 

\begin{figure}[ht]
	\centering
	\includegraphics[width=6.5truecm]{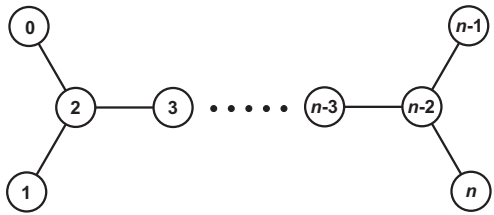}
	\caption{An illustration of the extended Dynkin diagram of type $\widetilde{D}_n$, $4 \leq n$, with its vertices labeled by integers from 0 to $n$.}
	\label{extDn}
\end{figure}

Finally, we give a brief description of {\it Dynkin diagrams} (see, e.\,g., \cite{hump}). These were introduced in the theory of Lie groups/algebras to describe particular sets of elements in lattices possessing integer quadratic forms -- so-called root systems. A Dynkin diagram is a graphical representation of the matrix of inner products of these roots. Given a root systems and its basis $S$, the vertices/nodes of its Dynkin diagram are the roots of $S$ and two nodes are not connected if the corresponding roots are orthogonal. If two nodes are not orthogonal, they are connected by one, two or three edges according as the angle between the corresponding roots is $2 \pi/3$, $3 \pi/4$ or $5 \pi/6$, respectively. In addition, a Dynkin diagram also encodes the lengths of roots. That is done
by marking the edge connecting two vertices whose corresponding roots are of different
length with an arrow pointing to the shorter root. Given a simple Lie algebra with its highest root $\gamma$, an extended root system is obtained by adding $-\gamma$ to the set of simple roots, which leads to the notion the {\it extended} Dynkin diagram of the algebra in question. In what follows, we will not use much of the properties of Dynkin diagrams as we will only be dealing with a particularly simple type of extended Dynkin diagrams, namely $\widetilde{D}_n$ ($4 \leq n$), depicted in Figure \ref{extDn}.

\section{Veldkamp Spaces of $\widetilde{D}_n$ and Two-Qubit Pauli Group}
$\widetilde{D}_n$, like any other graph, can be viewed/interpreted as a particular point-line incidence structure, 
$\mathcal{C}(\widetilde{D}_n)$, whose points and lines are, respectively, vertices and edges of $\widetilde{D}_n$; it thus features $n+1$ points and $n$ lines, where, for $5 \leq n$,  four points are of order one, two of order three and the remaining $n-5$ points being of order two. Let us adopt this view and have a detailed look at properties of the Veldkamp space of $\mathcal{C}(\widetilde{D}_n)$. We shall carry out this task step by step for $4 \leq n \leq 8$ in order to see how naturally the two-qubit (and, at the end, also three-qubit) Pauli group enters the stage.

\subsection{Case $n=4$}
\begin{figure}[ht]
	\centering
	\includegraphics[width=5.5truecm]{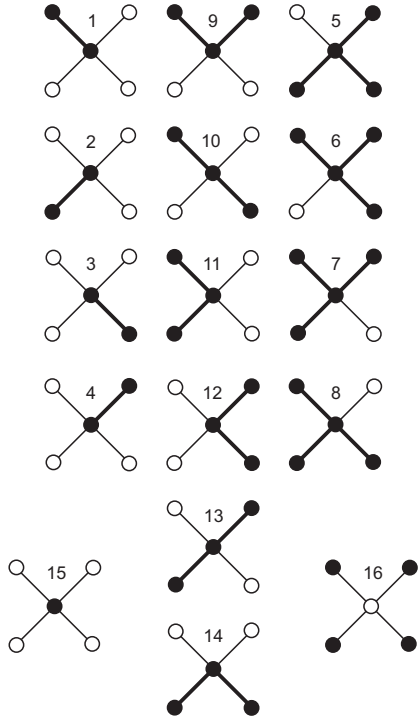}
	\caption{A diagrammatical representation of geometric hyperplanes of $\mathcal{C}(\widetilde{D}_4)$. Here, and in the sequel, a point of a hyperplane is represented by a filled circle and a line is drawn heavy if both of its points lie in a hyperplane.}
	\label{extD4-fig2}
\end{figure}

As readily discerned from Figure \ref{extDn} restricted to $n=4$, $\mathcal{C}(\widetilde{D}_4)$ consists of five points (0, 1, \dots, 4) and four lines, namely
$\{0,2\}$, $\{1,2\}$, $\{2,3\}$ and $\{2,4\}$, and features altogether 16 different geometric hyperplanes as listed in Table 1 and portrayed in Figure  \ref{extD4-fig2}. We see that each hyperplane except for the last one contains point 2; moreover, $H_{15}$ consists solely of this particular point and is thus contained in all preceding 14 hyperplanes. We further note that any other point of 
$\mathcal{C}(\widetilde{D}_4)$ is located in eight hyperplanes.

\begin{table}[h]
\begin{center}
\caption{The composition of 16 geometric hyperplanes $H_i$,  $1\leq i \leq 16$, of $\mathcal{C}(\widetilde{D}_4)$. The `+' symbol indicates which point of $\mathcal{C}(\widetilde{D}_4)$ (i.\,e., which vertex of $\widetilde{D}_4$) lies in a given hyperplane; for example, hyperplane $H_{10}$ consists of points 0, 2 and 4.} 
\vspace*{0.4cm} 
\resizebox{\columnwidth}{!}{%
\begin{tabular}{||c|cccccccccccccccc||}
\hline \hline
     &$H_1$&$H_2$&$H_3$&$H_4$&$H_5$&$H_6$&$H_7$&$H_8$&$H_9$&$H_{10}$&$H_{11}$&$H_{12}$&$H_{13}$&$H_{14}$&$H_{15}$&$H_{16}$  \\
\hline
~0~  &  +  &     &     &     &     &  +  &  +  &  +  &  +  &  +     &   +    &        &        &        &        &   +      \\
~1~  &     &  +  &     &     &  +  &     &  +  &  +  &     &        &   +    &        &   +    &   +    &        &   +      \\
~2~  &  +  &  +  &  +  &  +  &  +  &  +  &  +  &  +  &  +  &  +     &   +    &   +    &   +    &   +    &  +     &          \\
~3~  &     &     &     &  +  &  +  &  +  &  +  &     &  +  &        &        &   +    &   +    &        &        &   +      \\
~4~  &     &     &  +  &     &  +  &  +  &     &  +  &     &  +     &        &   +    &        &   +    &        &   +      \\
 \hline \hline
\end{tabular}%
}
\end{center}
\end{table} 

\begin{figure}[pth!]
	\centering
	\centerline{\includegraphics[width=5.5truecm]{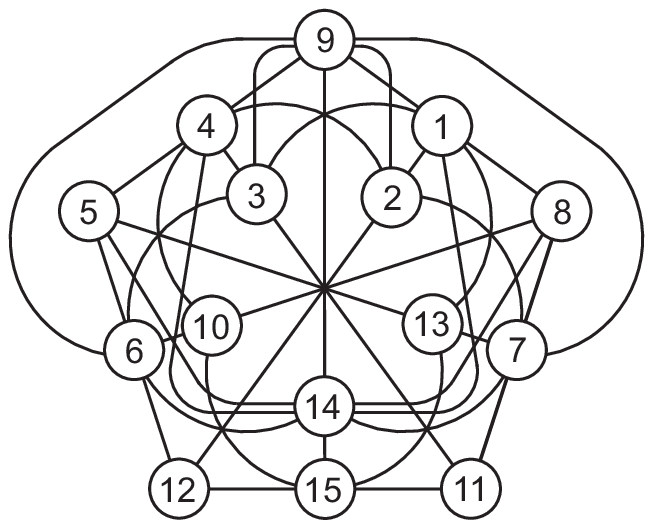}\hspace*{0.5cm}\includegraphics[width=5.5truecm]{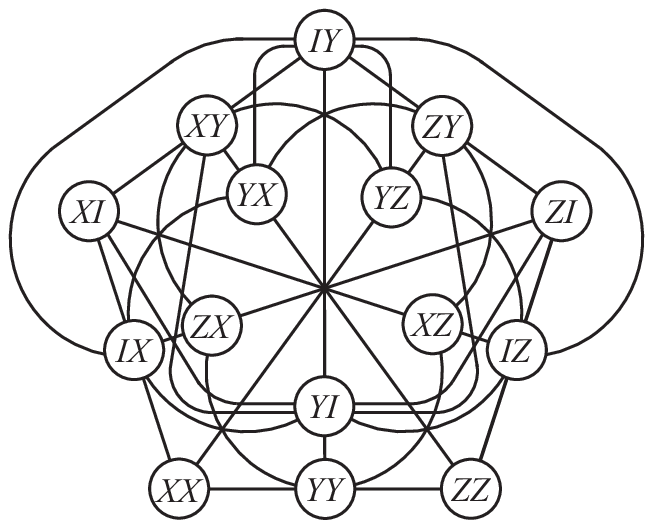}}
	\caption{A diagrammatical illustration of the fact that the point-line incidence structure encoded in Table 2 is isomorphic 
	to PG$(3,2)$ (left) and labeling the points of this space by the elements of the two-qubit Pauli group (right).}
	\label{extD4-fig3ab}
\end{figure}

From Table 1 (or Figure  \ref{extD4-fig2}) one infers that the Veldkamp space of $\mathcal{C}(\widetilde{D}_4)$ is endowed with 35 lines of size three, as given in Table 2. One sees that no such Veldkamp line contains $H_{16}$, the latter being thus regarded as an exceptional Veldkamp point of the geometry. It can readily be verified that the remaining 15 hyperplanes and all 35 Veldkamp lines form the projective space PG$(3,2)$. This is also illustrated in Figure  \ref{extD4-fig3ab}, {\it left}, employing, after Polster  \cite{pol}, 
a pictorial representation of PG$(3,2)$ built  around the pentagonal model (frequenly called the `doily') of the symplectic polar space $W(3,2)$ whose 15 lines are illustrated by straight-line-segments (10 of them) and/or arcs of circles (5). The remaining 20 lines of PG$(3,2)$ fall into four distinct orbits under the displayed automorphism of order five; from each orbit only one representative line is shown, namely $\{H_1, H_4, H_{14}\}$, $\{H_2, H_3, H_9\}$, $\{H_5, H_8, H_{14}\}$ and $\{H_6, H_7, H_9\}$, since the remaining ones can be readily obtained by rotating the figure through 72 degrees around the center of the pentagon.   

\begin{table}[t]
\begin{center}
\caption{The 35 three-point Veldkamp lines of $\mathcal{C}(\widetilde{D}_4)$. As in Table 1, the `+' symbol indicates which geometric hyperplane of $\mathcal{C}(\widetilde{D}_4)$ belongs to a given Veldkamp line.} 
\vspace*{0.4cm} 
\resizebox{\columnwidth}{!}{%
\begin{tabular}{||c|cccccccccccccccc||}
\hline \hline
     &$H_1$&$H_2$&$H_3$&$H_4$&$H_5$&$H_6$&$H_7$&$H_8$&$H_9$&$H_{10}$&$H_{11}$&$H_{12}$&$H_{13}$&$H_{14}$&$H_{15}$&$H_{16}$  \\
\hline
~1~  &  +  &  +  &     &     &     &     &     &     &     &        &        &    +   &        &        &        &         \\
~2~  &  +  &     &  +  &     &     &     &     &     &     &        &        &        &    +   &        &        &         \\
~3~  &  +  &     &     &  +  &     &     &     &     &     &        &        &        &        &    +   &        &         \\
~4~  &  +  &     &     &     &  +  &     &     &     &     &        &        &        &        &        &   +    &         \\
~5~  &  +  &     &     &     &     &  +  &     &     &     &        &   +    &        &        &        &        &         \\
~6~  &  +  &     &     &     &     &     &  +  &     &     &    +   &        &        &        &        &        &         \\
~7~  &  +  &     &     &     &     &     &     &  +  &  +  &        &        &        &        &        &        &         \\
~8~  &     &  +  &  +  &     &     &     &     &     &  +  &        &        &        &        &        &        &         \\
~9~  &     &  +  &     &  +  &     &     &     &     &     &    +   &        &        &        &        &        &         \\
~10~ &     &  +  &     &     &  +  &     &     &     &     &        &   +    &        &        &        &        &         \\
~11~ &     &  +  &     &     &     &  +  &     &     &     &        &        &        &        &        &    +   &         \\
~12~ &     &  +  &     &     &     &     &  +  &     &     &        &        &        &        &    +   &        &         \\
~13~ &     &  +  &     &     &     &     &     &  +  &     &        &        &        &    +   &        &        &         \\
~14~ &     &     &  +  &  +  &     &     &     &     &     &        &   +    &        &        &        &        &         \\
~15~ &     &     &  +  &     &  +  &     &     &     &     &   +    &        &        &        &        &        &         \\
~16~ &     &     &  +  &     &     &  +  &     &     &     &        &        &        &        &    +   &        &         \\
~17~ &     &     &  +  &     &     &     &  +  &     &     &        &        &        &        &        &    +   &         \\
~18~ &     &     &  +  &     &     &     &     &  +  &     &        &        &    +   &        &        &        &         \\
~19~ &     &     &     &  +  &  +  &     &     &     &  +  &        &        &        &        &        &        &         \\
~20~ &     &     &     &  +  &     &  +  &     &     &     &        &        &        &    +   &        &        &         \\
~21~ &     &     &     &  +  &     &     &  +  &     &     &        &        &    +   &        &        &        &         \\
~22~ &     &     &     &  +  &     &     &     &  +  &     &        &        &        &        &        &    +   &         \\
~23~ &     &     &     &     &  +  &  +  &     &     &     &        &        &    +   &        &        &        &         \\
~24~ &     &     &     &     &  +  &     &  +  &     &     &        &        &        &    +   &        &        &         \\
~25~ &     &     &     &     &  +  &     &     &  +  &     &        &        &        &        &    +   &        &         \\
~26~ &     &     &     &     &     &  +  &  +  &     &  +  &        &        &        &        &        &        &         \\
~27~ &     &     &     &     &     &  +  &     &  +  &     &    +   &        &        &        &        &        &         \\
~28~ &     &     &     &     &     &     &  +  &  +  &     &        &    +   &        &        &        &        &         \\
~29~ &     &     &     &     &     &     &     &     &  +  &    +   &    +   &        &        &        &        &         \\
~30~ &     &     &     &     &     &     &     &     &  +  &        &        &    +   &    +   &        &        &         \\
~31~ &     &     &     &     &     &     &     &     &  +  &        &        &        &        &    +   &    +   &         \\
~32~ &     &     &     &     &     &     &     &     &     &    +   &        &    +   &        &    +   &        &         \\
~33~ &     &     &     &     &     &     &     &     &     &    +   &        &        &    +   &        &    +   &         \\
~34~ &     &     &     &     &     &     &     &     &     &        &    +   &    +   &        &        &    +   &         \\
~35~ &     &     &     &     &     &     &     &     &     &        &    +   &        &    +   &    +   &        &         \\
\hline \hline
\end{tabular}%
}
\end{center}
\end{table} 
\begin{table}[pth!]
\begin{center}
\caption{The 15 hyperplanes of $\mathcal{C}(\widetilde{D}_4)$ are in a one-to-one correspondence with the 15 elements of $\overline{{\cal P}}_2$; note that the `exceptional' hyperplane $H_{16}$ corresponds to the same group element as  $H_{15}$, because the two hyperlanes are complementary.} 
\vspace*{0.4cm} 
\resizebox{\columnwidth}{!}{%
\begin{tabular}{||c|ccccccccccccccc|c||}
\hline \hline
$\mathcal{C}(\widetilde{D}_4)$     
     &$H_1$&$H_2$&$H_3$&$H_4$&$H_5$&$H_6$&$H_7$&$H_8$&$H_9$&$H_{10}$&$H_{11}$&$H_{12}$&$H_{13}$&$H_{14}$&$H_{15}$&$H_{16}$  \\
\hline
$\overline{{\cal P}}_2$     
     & $ZY$&$YZ$ &$YX$ &$XY$ &$XI$ &$IX$ &$IZ$ &$ZI$ &$IY$ &$ZX$    &$ZZ$    &$XX$    &$XZ$    &$YI$    &$YY$    &$YY$      \\
 \hline \hline
\end{tabular}%
}
\end{center}
\end{table} 

Let us now return back to our $\widetilde{D}_4$ and label its five vertices by five distinct elements of the two-qubit Pauli group, 
$\overline{{\cal P}}_2$. One can take any five elements requiring only that their product equals $II$; this constraint is necessary to ensure that the induced labeling of the points the associated Veldkamp space has the property that the product of any three collinear elements is also equal to $II$. Given the symmetry of $\widetilde{D}_4$, one of the most natural choices is 
 ($A \otimes B$ is shorthanded to $AB$ in the sequel) 
$$0 \rightarrow XI,~ 1 \rightarrow IX,~ 2 \rightarrow YY,~ 3 \rightarrow ZI,~ 4 \rightarrow IZ. $$
Assume further that each hyperplane acquires the label that is the product of the group elements attached to the points it consists of; thus, for example, $H_1$, comprising points 0 and 2, will bear the label $(XI).(YY) = ZY$. Hence, we arrive at a {\it one-to-one} correspondence between the 15 Veldkamp points of $\mathcal{C}(\widetilde{D}_4)$ and the 15 elements of $\overline{{\cal P}}_2$ as shown in Table 3. Even more interestingly, employing Table 2 one can check that not only is the product of three elements on each Veldkamp line equal to the identity element of the group, but -- as handy rendered by Figure \ref{extD4-fig3ab}, {\it right} -- the three elements that lie on a line of the selected copy\footnote{There are altogether 28 distinct copies of $W(3,2)$ contained in PG$(3,2)$.} of $W(3,2)$ our PG$(3,2)$ was built around pairwise commute.

\begin{figure}[t]
	\centering
	\includegraphics[width=12.0truecm]{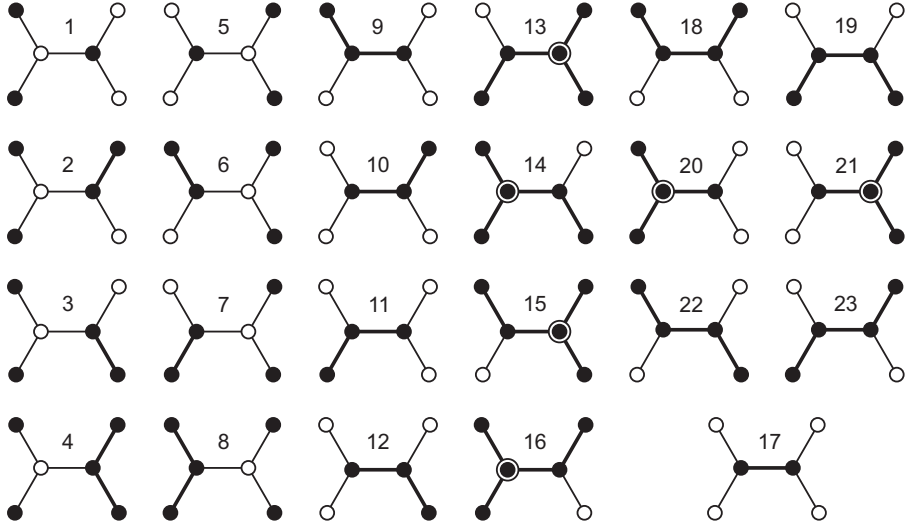}
	\caption{Geometric hyperplanes of $\mathcal{C}(\widetilde{D}_5)$. }
	\label{extD5-fig2}
\end{figure}

\begin{table}[h]
\begin{center}
\caption{A `double-six' of Veldkamp lines generated by hyperplanes of the first family.} 
\vspace*{0.4cm} 
\resizebox{\columnwidth}{!}{%
\begin{tabular}{||c|cccc|cccc|cccccc||}
\hline \hline
     &$H_1$&$H_2$&$H_3$&$H_4$&$H_5$&$H_6$&$H_7$&$H_8$&$H_{13}$&$H_{14}$&$H_{15}$&$H_{16}$&$H_{20}$&$H_{21}$  \\
\hline
~1~  &  +  &  +  &     &     &     &     &     &     &        &    +   &        &        &        &          \\
~2~  &  +  &     &  +  &     &     &     &     &     &        &        &        &    +   &        &          \\
~3~  &  +  &     &     &  +  &     &     &     &     &        &        &        &        &    +   &          \\
~4~  &     &  +  &  +  &     &     &     &     &     &        &        &        &        &    +   &          \\
~5~  &     &  +  &     &  +  &     &     &     &     &        &        &        &    +   &        &          \\
~6~  &     &     &  +  &  +  &     &     &     &     &        &    +   &        &        &        &          \\
\hline
~7~  &     &     &     &     &  +  &  +  &     &     &   +    &        &        &        &        &          \\
~8~  &     &     &     &     &  +  &     &  +  &     &        &        &   +    &        &        &          \\
~9~  &     &     &     &     &  +  &     &     &  +  &        &        &        &        &        &  +       \\
~10~ &     &     &     &     &     &  +  &  +  &     &        &        &        &        &        &  +       \\
~11~ &     &     &     &     &     &  +  &     &  +  &        &        &   +    &        &        &          \\
~12~ &     &     &     &     &     &     &  +  &  +  &   +    &        &        &        &        &          \\
\hline \hline
\end{tabular}%
}
\end{center}
\end{table}

\subsection{Case $n=5$}
\begin{figure}[h]
	\centering
	\includegraphics[width=10.0truecm]{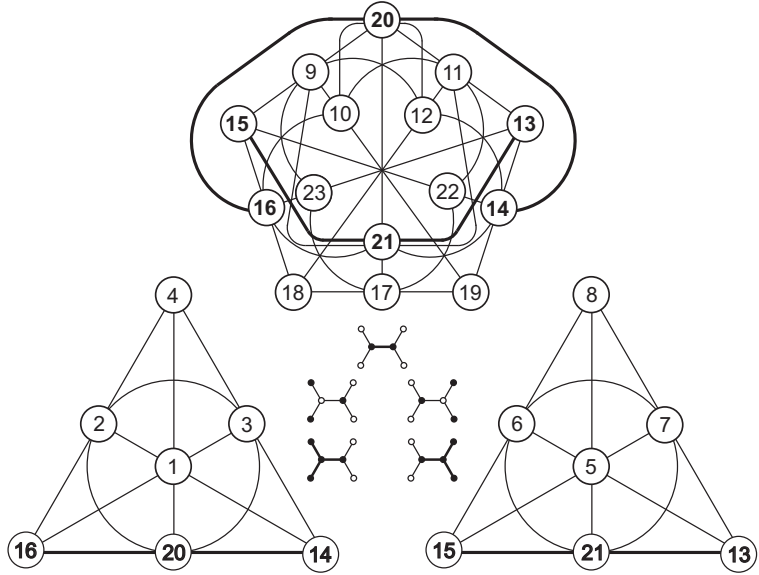}
	\caption{A diagrammatic illustration of all projective subgeometries of the Veldkamp space of $\mathcal{C}(\widetilde{D}_5)$. The central inset depicts those geometric hyperplanes each of which fully encodes all Veldkamp points of a particular subgeometry, namely of the PG$(3,2)$ ($H_{17}$, top), the two Fano planes ($H_1$ and $H_5$, middle) as well as of the two shared lines ($H_{20}$ and $H_{21}$, bottom).}
	\label{extD5-fig4}
\end{figure}

$\mathcal{C}(\widetilde{D}_5)$ contains six points (0, 1, \dots, 5) and five lines ($\{0,2\}$, $\{1,2\}$, $\{2,3\}$, 
$\{3,4\}$, and $\{3,5\}$). Its 23 geometric hyperplanes, shown in Figure \ref{extD5-fig2},  can be split into two disjoint families, namely $\{H_1, H_2, \dots, H_8\}$ and $\{H_9, H_{10}, \dots, H_{23}\}$ according as they do not or do contain the line $\{2,3\}$, respectively. The former family can further be divided into two subfamilies, $\{H_1, H_2, H_3, H_4\}$ and $\{H_5, H_6, H_7, H_8\}$, depending on whether a hyperplane misses, respectively, point 2 or point 3.
These splittings have a deep geometrical meaning once we see all 47 three-point Veldkamp lines $\mathcal{C}(\widetilde{D}_5)$ is found to possess. Twelve of them are generated by hyperplanes of the first family and they are given in Table 4; the remaining 35 are defined by hyperplanes of the second family and their properties are summarized in Table 5.
As it is obvious from Table 4, the four hyperplanes of either subfamily define PG$(2,2)$, the Fano plane, with one line omitted; the latter being line $\{H_{14}, H_{16}, H_{20}\}$ for the first and $\{H_{13}, H_{15}, H_{21}\}$ for the second subfamily. Comparing Table 4 with Figure \ref{extD5-fig2} we see that the Veldkamp points of the first/second Fano plane are those seven hyperplanes that contain 
$H_1$/$H_5$, and the two `missing' lines consist of those three hyperplanes that contain $H_{20}$/$H_{21}$.  The 15 Veldkamp points (that are all geometric hyperplanes incorporating $H_{17}$, see Figure \ref{extD5-fig2}) and 35 Veldkamp lines of Table 5 define a projective space isomorphic to PG$(3,2)$; this space also contains two `missing' lines (marked in italics) of the above-described Fano planes. Summing up, the subgeometry of the Veldkamp space of $\mathcal{C}(\widetilde{D}_5)$ with lines of size three comprises the projective space PG$(3,2)$ and a couple of disjoint Fano planes, each sharing a line with this PG$(3,2)$ -- as depicted in Figure \ref{extD5-fig4}. (It is also worth adding that in this case we have no exceptional Veldkamp point(s) since there is no geometric hyperplane lacking line 
$\{2,3\}$.)

\begin{table}[t]
\begin{center}
\caption{The 35 Veldkamp lines generated by 15 hyperplanes of the second family.} 
\vspace*{0.4cm} 
\resizebox{\columnwidth}{!}{%
\begin{tabular}{||c|ccccccccccccccc||}
\hline \hline
     &$H_9$&$H_{10}$&$H_{11}$&$H_{12}$&$H_{13}$&$H_{14}$&$H_{15}$&$H_{16}$&$H_{17}$&$H_{18}$&$H_{19}$&$H_{20}$&$H_{21}$&$H_{22}$&$H_{23}$  \\
\hline
~13~ &  +  &  +     &        &        &        &        &        &        &        &        &   +    &        &        &        &     \\
~14~ &  +  &        &  +     &        &        &        &        &        &        &        &        &        &    +   &        &      \\
~15~ &  +  &        &        &  +     &        &        &        &        &        &        &        &        &        &        & +    \\
~16~ &  +  &        &        &        &  +     &        &        &        &    +   &        &        &        &        &        &      \\
~17~ &  +  &        &        &        &        &  +     &        &        &        &   +    &        &        &        &        &     \\
~18~ &  +  &        &        &        &        &        &  +     &        &        &        &        &   +    &        &        &      \\
~19~ &  +  &        &        &        &        &        &        &  +     &        &        &        &        &        &    +   &      \\
~20~ &     &  +     &  +     &        &        &        &        &        &        &        &        &        &        &    +   &      \\
~21~ &     &  +     &        &  +     &        &        &        &        &        &        &        &   +    &        &        &      \\
~22~ &     &  +     &        &        &  +     &        &        &        &        &   +    &        &        &        &        &       \\
~23~ &     &  +     &        &        &        &  +     &        &        &   +    &        &        &        &        &        &      \\
~24~ &     &  +     &        &        &        &        &  +     &        &        &        &        &        &        &        &  +   \\
~25~ &     &  +     &        &        &        &        &        &  +     &        &        &        &        &    +   &        &      \\
~26~ &     &        &  +     &  +     &        &        &        &        &        &   +    &        &        &        &        &      \\
~27~ &     &        &  +     &        &  +     &        &        &        &        &        &        &    +   &        &        &      \\
~28~ &     &        &  +     &        &        &  +     &        &        &        &        &        &        &        &        & +    \\
~29~ &     &        &  +     &        &        &        &  +     &        &   +    &        &        &        &        &        &     \\
~30~ &     &        &  +     &        &        &        &        &  +     &        &        &   +    &        &        &        &     \\
~31~ &     &        &        &  +     &  +     &        &        &        &        &        &        &        &        &   +    &     \\
~32~ &     &        &        &  +     &        &  +     &        &        &        &        &        &        &    +   &        &     \\
~33~ &     &        &        &  +     &        &        &  +     &        &        &        &   +    &        &        &        &     \\
~34~ &     &        &        &  +     &        &        &        &  +     &   +    &        &        &        &        &        &     \\
~35~ &     &        &        &        &  +     &  +     &        &        &        &        &   +    &        &        &        &     \\
~{\it 36}~ &     &        &        &        &  +     &        &  +     &        &        &        &        &        &   +    &        &     \\
~37~ &     &        &        &        &  +     &        &        &  +     &        &        &        &        &        &        & +   \\
~38~ &     &        &        &        &        &  +     &  +     &        &        &        &        &        &        &    +   &     \\
~{\it 39}~ &     &        &        &        &        &  +     &        &  +     &        &        &        &    +   &        &        &     \\
~40~ &     &        &        &        &        &        &  +     &  +     &        &    +   &        &        &        &        &     \\
~41~ &     &        &        &        &        &        &        &        &  +     &    +   &    +   &        &        &        &     \\
~42~ &     &        &        &        &        &        &        &        &  +     &        &        &    +   &    +   &        &     \\
~43~ &     &        &        &        &        &        &        &        &  +     &        &        &        &        &    +   &    +\\
~44~ &     &        &        &        &        &        &        &        &        &    +   &        &    +   &        &    +   &     \\
~45~ &     &        &        &        &        &        &        &        &        &    +   &        &        &    +   &        &    +\\
~46~ &     &        &        &        &        &        &        &        &        &        &    +   &    +   &        &        &    + \\
~47~ &     &        &        &        &        &        &        &        &        &        &    +   &        &    +   &    +   &      \\
\hline \hline
\end{tabular}%
}
\end{center}
\end{table} 

\begin{table}[h]
\begin{center}
\caption{Labeling the Veldkamp points of $\mathcal{C}(\widetilde{D}_5)$ by the elements of the two-qubit Pauli group.} 
\vspace*{0.4cm} 
\resizebox{\columnwidth}{!}{%
\begin{tabular}{||l|ccccccccccccccc||}
\hline \hline
$\overline{{\cal P}}_2$
         & $XY$& $YX$   & $ZY$   & $YZ$   &  $ZI$  &  $IZ$  &  $XI$  &  $IX$  & $YY$   & $XX$   & $ZZ$   &$IY$    & $YI$   & $XZ$   & $ZX$ \\
				\hline
PG$(3,2)$&$H_9$&$H_{10}$&$H_{11}$&$H_{12}$&$H_{13}$&$H_{14}$&$H_{15}$&$H_{16}$&$H_{17}$&$H_{18}$&$H_{19}$&$H_{20}$&$H_{21}$&$H_{22}$&$H_{23}$  \\
\hline
1st Fano &     &$H_2$   &        &$H_3$   &        &        &        &        &$H_1$   &        &        &        &$H_4$   &        &  \\
2nd Fano &$H_6$&        &$H_7$   &        &        &        &        &        &$H_5$   &        &        &$H_8$   &        &        &  \\
\hline \hline
\end{tabular}%
}
\end{center}
\end{table} 

\begin{figure}[pth!]
	\centering
	\includegraphics[width=10.0truecm]{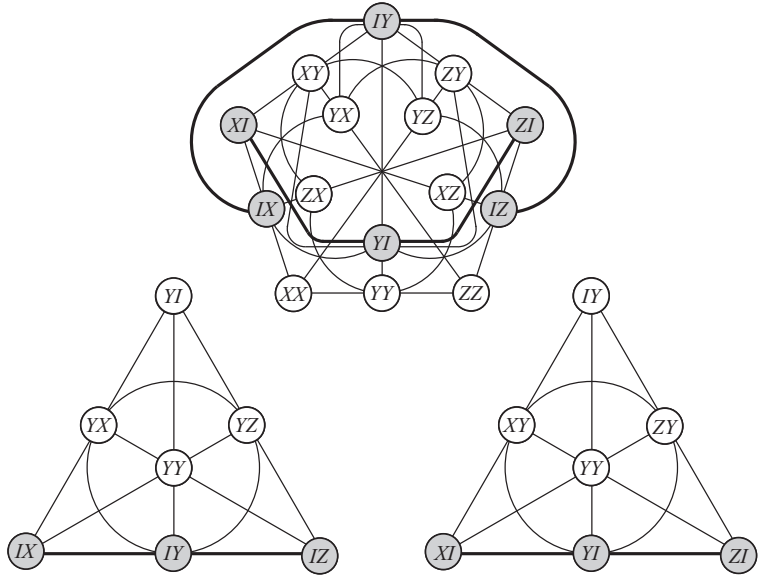}
	\caption{A general view of the projective subgeometries of the Veldkamp space of $\mathcal{C}(\widetilde{D}_5)$ in terms of the elements of the two-qubit Pauli group.}
	\label{extD5-fig4b}
\end{figure}

Next, pursuing the strategy of the preceding case, one labels six vertices of $\widetilde{D}_5$ by elements of the two-qubit Pauli group as
$$0 \rightarrow ZI,~ 1 \rightarrow XI,~ 2 \rightarrow YI,~ 3 \rightarrow IY,~ 4 \rightarrow IZ,~ 5 \rightarrow IX, $$
and, making use of Figure \ref{extD5-fig2}, gets a bijection between the elements of the group and the points of the PG$(3,2)$  in the form shown in Table 6; this table also shows which elements of the two-qubit Pauli group are ascribed to the remaining eight points of the two Fano planes. Figure \ref{extD5-fig4b} renders a pictorial representation of the information gathered in Table 6. 
The bijection is of similar nature as the one of the $n=4$ case: that is, a line of the PG$(3,2)$ entails three group elements whose product is $II$, and a line of the distinguished copy of $W(3,2)$ gathers a triple of mutually commuting elements.
Yet, it also features an interesting novelty due to the fact that our PG$(3,2)$ has two distinguished lines that it shares with the two Fano planes. If we forget about the six group elements located on these two lines (highlighted by light shading in Figure \ref{extD5-fig4b}) we shall find that the remaining nine elements form within the $W(3,2)$ nothing but a $3 \times 3$ grid (or, what amounts to the same, a copy of the generalized quadrangle GQ$(2,1)$). Physical importance of this observation lies with the fact \cite{ps} that any such grid with the labeling inherited from that of its parent $W(3,2)$ represents the so-called Mermin-Peres magic square -- one of the simplest proofs of the  (Bell-)Kochen-Specker theorem first proposed by Mermin \cite{mer} and Peres \cite{per}.

\begin{figure}[pth!]
	\centering
	\includegraphics[width=12.0truecm]{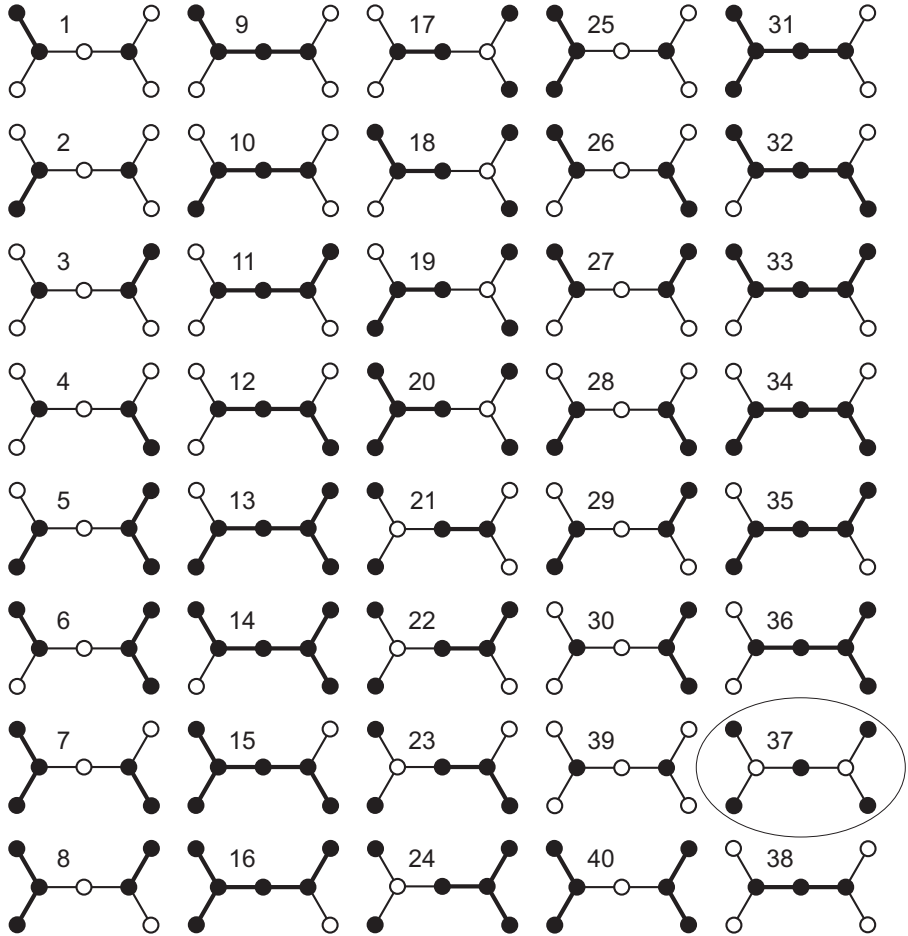}
	\caption{Geometric hyperplanes of $\mathcal{C}(\widetilde{D}_6)$; an ellipse marks the (single) exceptional hyperplane.}
	\label{extD6-fig2}
\end{figure}

\subsection{Case $n=6$}
\begin{figure}[t]
	\centering
	\includegraphics[width=5.0truecm]{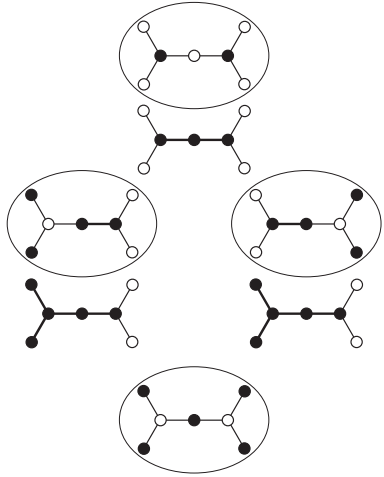}
	\caption{A symbolic structure of the Veldkamp space of $\mathcal{C}(\widetilde{D}_6)$. Each projective space (starting with the 
	PG$(4,2)$ at the top and ending with the `exceptional' PG$(1,2)$ at the bottom) is represented by a single hyperplane, viz. the one that fully determines the remaining hyperplanes defining the space in question. Marked by ellipses are those spaces that are not properly contained in any other space.}
	\label{extD6-fig3}
\end{figure}

\begin{figure}[pth!]
	\centering
	\centerline{\includegraphics[width=6.0truecm]{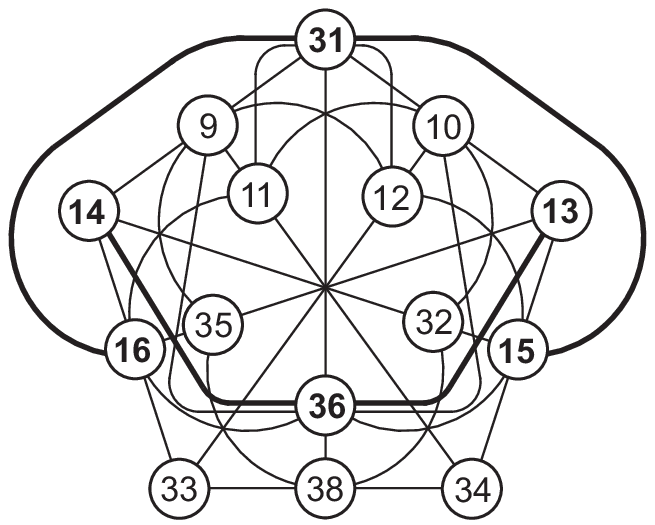}\hspace*{0.5cm}\includegraphics[width=6.0truecm]{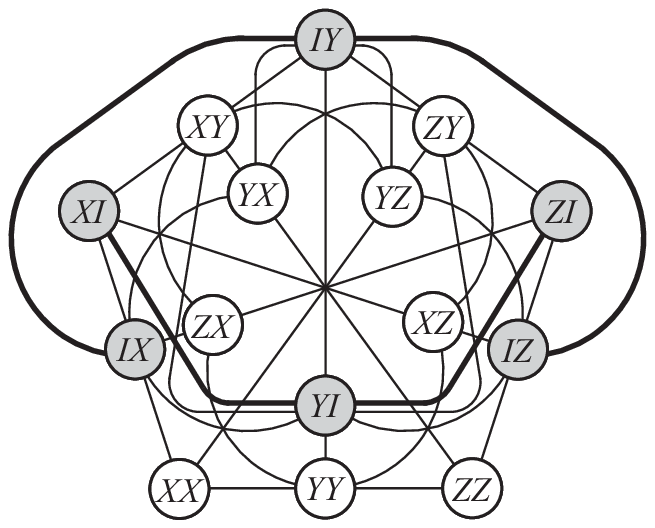}}
	\caption{{\it Left:} The structure of the distinguished PG$(3,2)$; the two lines shown in bold are those shared with the two Fano planes, one in each. {\it Right:} The two different two-qubit labelings of the vertices of  $\widetilde{D}_6$ lead to the identical labelings of the points of  this space. As the labeling is the same as that of the $n=5$ case (see Figure \ref{extD5-fig4b}), removal of the highlighted elements leads to the same Mermin-Peres magic square as in the preceding case.}
	\label{extD6-fig4}
\end{figure}

\begin{figure}[pth!]
	\centering
	\includegraphics[width=12.0truecm]{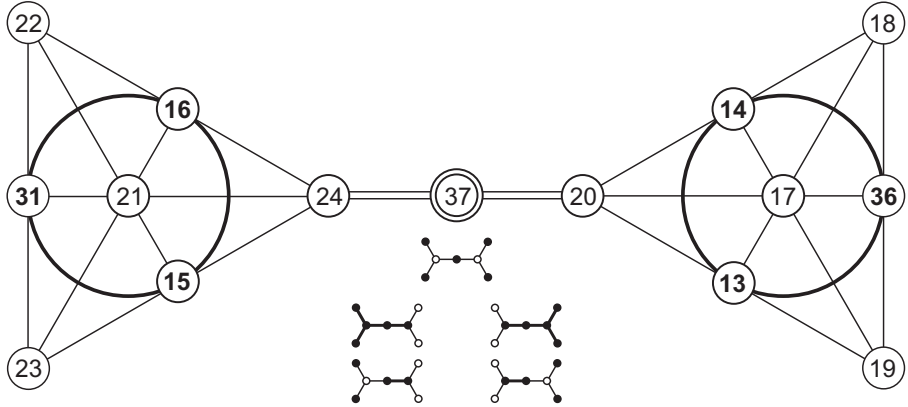}
	\caption{The two Fano planes and the exceptional Veldkamp line interconnecting them. The lines that also belong to the distinguished PG$(3,2)$ are drawn in boldface. (As in Figure \ref{extD5-fig4}, the inset shows the representative hyperplanes for these projective geometries). }
	\label{extD6-fig5}
\end{figure}

$\mathcal{C}(\widetilde{D}_6)$, comprising seven points ($0, 1, \dots, 6$) and six lines ($\{0,2\}$, $\{1,2\}$, $\{2,3\}$, 
$\{3,4\}$, $\{4,5\}$, and $\{4,6\}$), is found to possess 40 geometric hyperplanes -- all depicted in Figure \ref{extD6-fig2} -- and as many as 168 Veldkamp lines of size three. We shall not give an explicit list of the latter here, since Figure \ref{extD6-fig2} contains all the information we need to find all projective subgeometries of the corresponding Veldkamp space. A detailed inspection of this figure leads to the following observations. The smallest hyperplane, $H_{39}$, is contained in other 30 hyperplanes, which together yield 155 three-point Veldkamp lines and form the projective space isomorphic to PG$(4,2)$ (see Sect. 2); this projective space also contains a distinguished copy of PG$(3,2)$ that is defined by $H_{38}$ and the other 14 hyperplanes containing it. Then we have a pair of complementary Fano planes, one defined by seven hyperplanes containing $H_{17}$, the other by seven hyperplanes comprising $H_{21}$. Either of the two Fano planes shares a line with the distinguished copy of PG$(3,2)$ (and, hence, also with the parent PG$(4,2)$); it is the line defined by three hyperplanes containing $H_{36}$ in the former and $H_{31}$ in the latter case. This already accounts for $155 + (2 \times 6) = 167$ Veldkamp lines.  The remaining Veldkamp line is $\{H_{20},H_{24},H_{37}\}$, the joint of the two Fano planes, which is the {\it only} size-three Veldkamp line passing through the exceptional Veldkamp point $H_{37}$. This hierarchy of projective spaces living within our Veldkamp space can be expressed in a succinct form as displayed in Figure \ref{extD6-fig3}, with its `simpler' parts being shown in full in Figure \ref{extD6-fig4}, {\it left},  and Figure \ref{extD6-fig5}.

\begin{figure}[pth!]
	\centering
	\centerline{\includegraphics[width=5.0truecm]{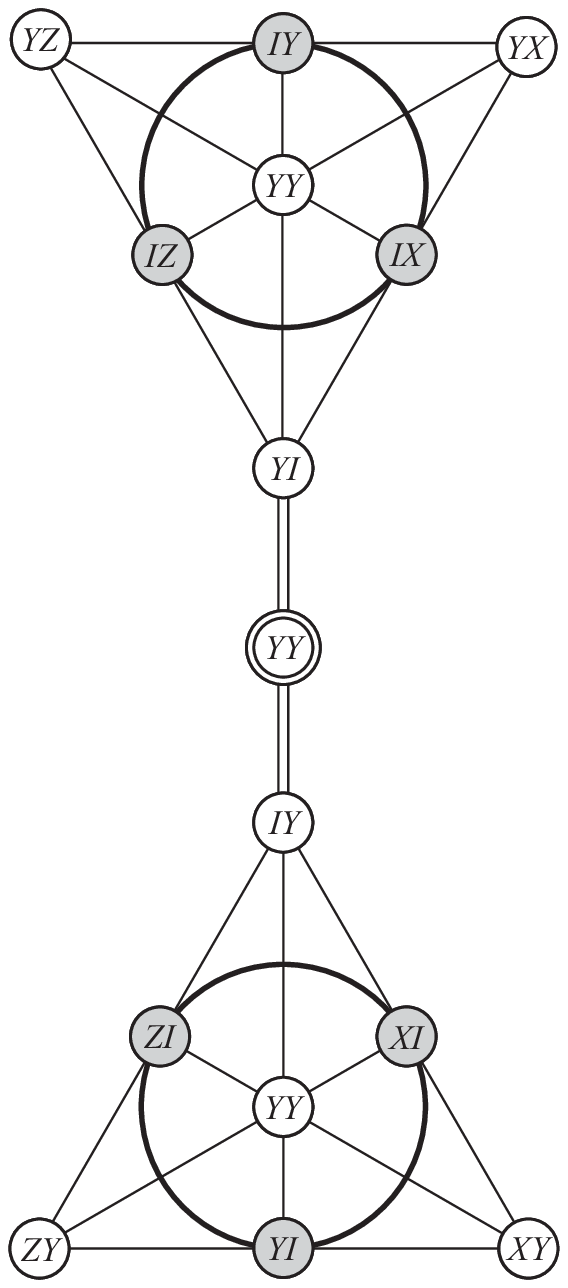}\hspace*{.5cm}
	\includegraphics[width=5.0truecm]{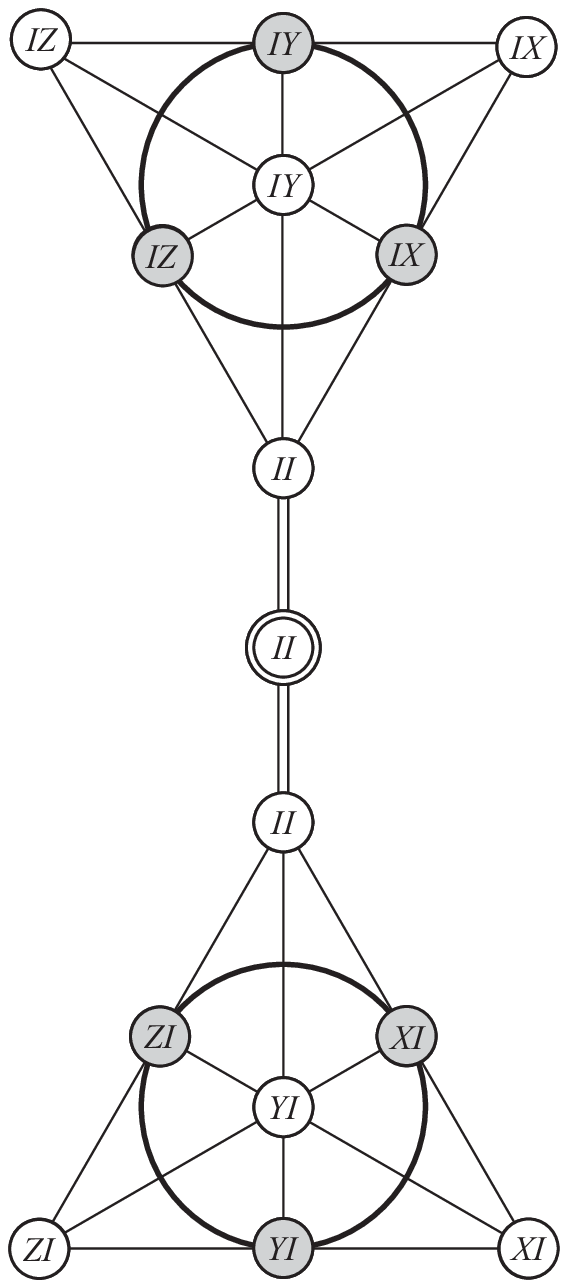}}
	\caption{The pair of interconnected Fano planes in light of the two distinct two-qubit labelings. }
	\label{extD6-fig5bc}
\end{figure}

An interesting thing here happens when it comes to the relation with the two-qubit Pauli group, as we have now at our disposal two natural labelings of the vertices of  $\widetilde{D}_6$, both featuring, unlike the previous two cases, also the identity element; 
in particular,
$$0 \rightarrow ZI,~ 1 \rightarrow XI,~ 2 \rightarrow YI,~ 3 \rightarrow II,~ 4 \rightarrow IY,~ 5 \rightarrow IZ,~ 6 \rightarrow IX, $$
and
$$0 \rightarrow ZI,~ 1 \rightarrow XI,~ 2 \rightarrow II,~ 3 \rightarrow YY,~ 4 \rightarrow II,~ 5 \rightarrow IZ,~ 6 \rightarrow IX. $$
Although these two labelings give one and the same labeling of the points of the distinguished PG$(3,2)$, that portrayed in Figure \ref{extD6-fig4}, {\it right}, there is a substantial difference between them when the two interconnected Fano planes are concerned, as readily discerned from  Figure \ref{extD6-fig5bc}; this difference is most pronounced for the exceptional Veldkamp line, as in the second case its three points acquire the same label, this being the identity element at that.

\begin{figure}[pth!]
	\centering
	\includegraphics[width=8.0truecm]{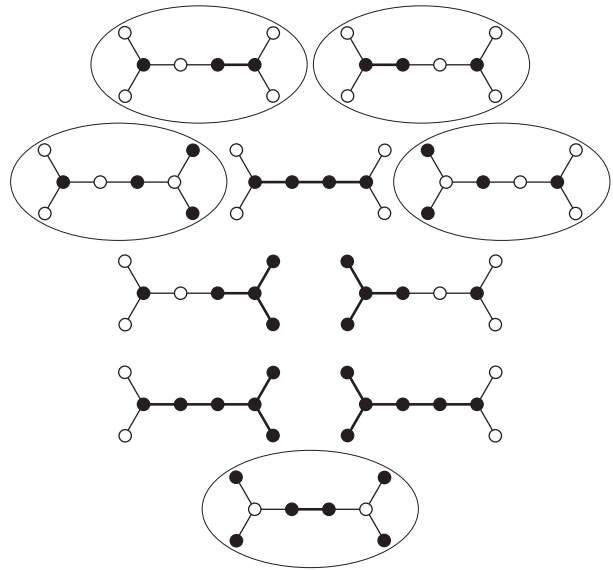}
	\caption{Stratification of the Veldkamp space of $\mathcal{C}(\widetilde{D}_7)$ in terms of projective spaces it contains.}
	\label{extD7-fig2}
\end{figure}

\subsection{Case $n=7$}
\begin{figure}[t]
	\centering
	\includegraphics[width=12.0truecm]{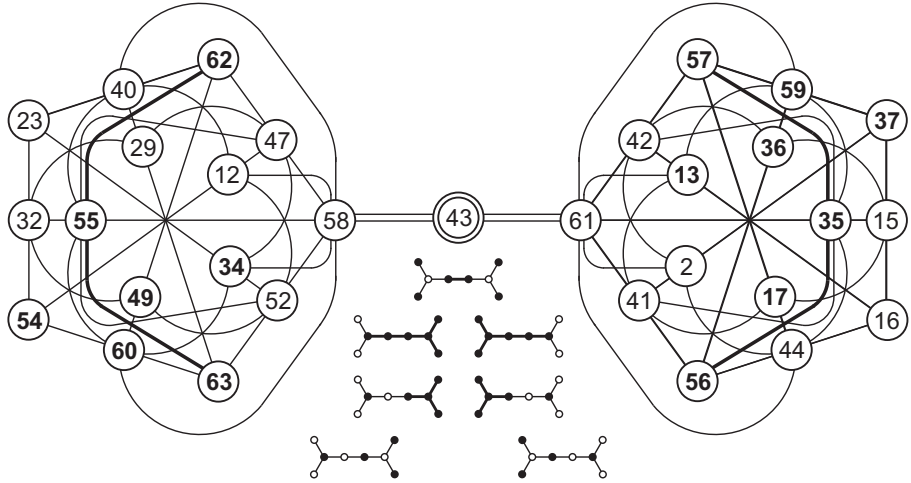}
	\caption{The left  and right  PG$(3,2)$s connected by the exceptional Veldkamp line. The seven points in either of these spaces that are numbered in boldface represent the Fano plane lying also in the corresponding PG$(4,2)$; the line of this Fano plane that also belongs to the distinguished PG(3,2) is drawn thick. (As before, the inset depicts the representative hyperplane for each of the projective spaces mentioned.)}
	\label{extD7-fig7}
\end{figure}

Our description of this case will be rather brief, as there are no conceptual novelties when compared with the cases we have already addressed.
Eight points ($0, 1, \dots, 7$) and seven lines ($\{0,2\}$, $\{1,2\}$, $\{2,3\}$, $\{3,4\}$, $\{4,5\}$, $\{5,6\}$ and $\{5,7\}$) of  
$\mathcal{C}(\widetilde{D}_7)$ are found to accommodate 64 geometric hyperplanes and 332 Veldkamp lines of size three, whose hierarchical projective structure is shown in Figure  \ref{extD7-fig2}. 
As in the preceding subsection (see Figure \ref{extD6-fig3}) each hyperplane of Figure \ref{extD7-fig2} represents a PG$(d,2)$ whose dimension $d$ is one less than is the number of points that are not contained in the hyperplane (i.e., the number of empty circles).
Thus, from the 1st row of the figure one reads off that maximum projective (sub)spaces of the Veldkamp space of $\mathcal{C}(\widetilde{D}_6)$ are two PG$(4,2)$s, referred simply as the left PG$(4,2)$ and the right PG$(4,2)$, which intersect in the distinguished PG$(3,2)$ (the 2nd row, the middle hyperplane). There are other two PG(3,2)s (let us also call them left and right) present here (the 2nd row, left and right hyperplanes). They are disjoint, but either of them shares a Fano plane with its PG(4,2) counterpart (namely the Fano plane  represented by the respective hyperplane shown in the 3rd row), and a line with the distinguished PG$(3,2)$ (the lines represented by the respective hyperplane shown in the 4th row). Finally, there is an `exceptional' projective line that is, similarly to the $n=6$ case, the only Veldkamp line passing through the sole exceptional hyperplane (shown at the bottom of the figure). The analogy with the previous case is even deeper, since this exceptional Veldkamp line also joins the left PG$(3,2)$ with the right one -- as shown in Figure \ref{extD7-fig7} (compare with Figure \ref{extD6-fig5}). To decipher the content of the last figure, one has to say a few words on how the numbering of hyperplanes was done. 
A hyperplane with a smaller number of points precedes a hyperplane with a larger number of points. If two hyperplanes have the same number of points, first goes that which contains a point labeled by the smallest integer. If also this point is the same in both hyperplanes, then first goes that featuring a point labeled by the second smallest integer, and so on. Here are shown a few hyperplanes\footnote{Note that, unlike the previous cases, our numbering starts here with `0'.}
$$H_0 = \{2,3,5\},$$  
$$H_1 = \{2,4,5\},$$ 
$$H_2 = \{0,1,3,5\},$$ 
$$H_3 = \{0,2,3,5\},$$ 
$$H_4 = \{0,2,4,5\},$$
$$\vdots~~\vdots~~\vdots$$
$$H_{43} = \{0,1,3,4,6,7\},$$
$$\vdots~~\vdots~~\vdots$$
$$H_{62} = \{0,2,3,4,5,6,7\},$$
$$H_{63} = \{1,2,3,4,5,6,7\}.$$

\begin{figure}[t]
	\centering
	\centerline{\includegraphics[width=6.5truecm]{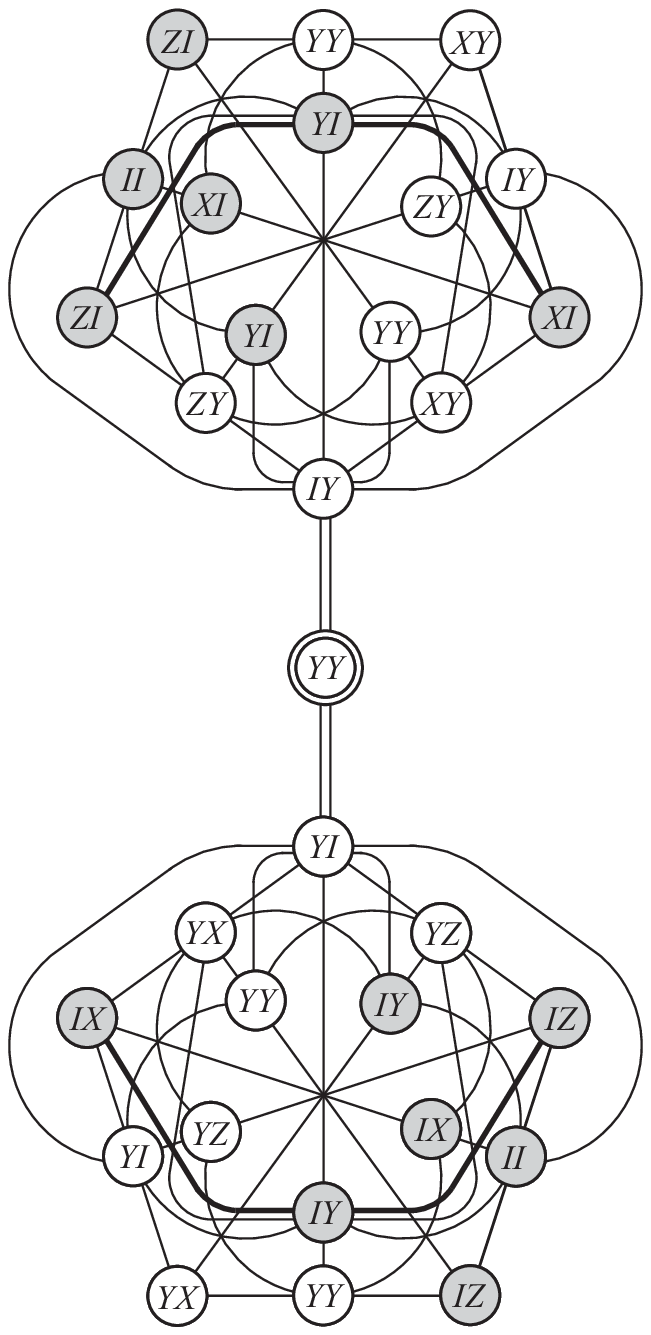}\hspace*{0.5cm}\includegraphics[width=6.5truecm]{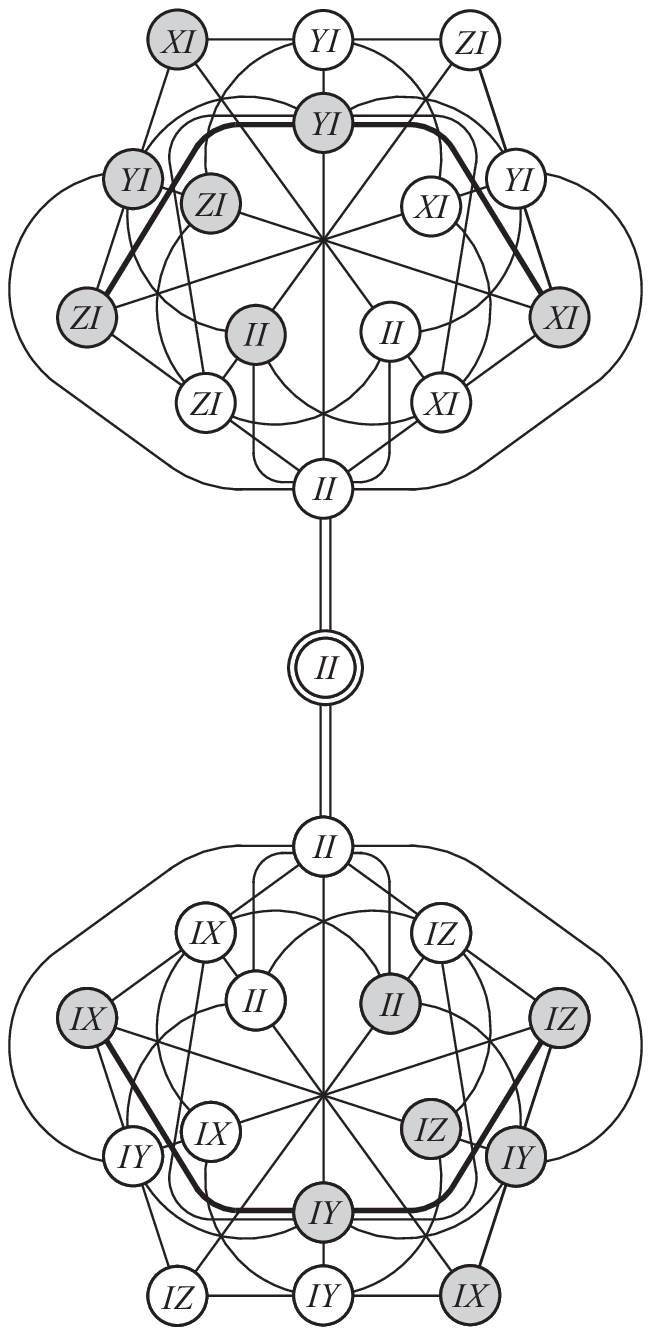}}
	\caption{The interconnected left  and right  PG$(3,2)$s in terms of the two two-qubit labelings.}
	\label{extD7-fig7bc}
\end{figure}

Regarding a relation with the two-qubit Pauli group, we again see close parallels with the $n=6$ case. For not only do we have again two natural labellings of the vertices of  $\widetilde{D}_7$,
$$0 \rightarrow ZI,~ 1 \rightarrow XI,~ 2 \rightarrow YI,~ 3 \rightarrow II,~ 4 \rightarrow II,~ 5 \rightarrow IY,~ 6 \rightarrow IZ,
~7 \rightarrow IX, $$
and
$$0 \rightarrow ZI,~ 1 \rightarrow XI,~ 2 \rightarrow II,~ 3 \rightarrow YI,~ 4 \rightarrow IY,~ 5 \rightarrow II,~ 6 \rightarrow IZ,
~7 \rightarrow IX, $$
but these also give identical labellings of the distinguished PG$(3,2)$, furnishing the same prominent copy of the Mermin-Peres magic square; in addition, as illustrated in Figure \ref{extD7-fig7bc}, the two labelings of the exceptional Veldkamp line are the same as those of its $n=6$ counterpart.

\subsection{Case $n=8$}
\begin{figure}[pth!]
	\centering
	\includegraphics[width=10.0truecm]{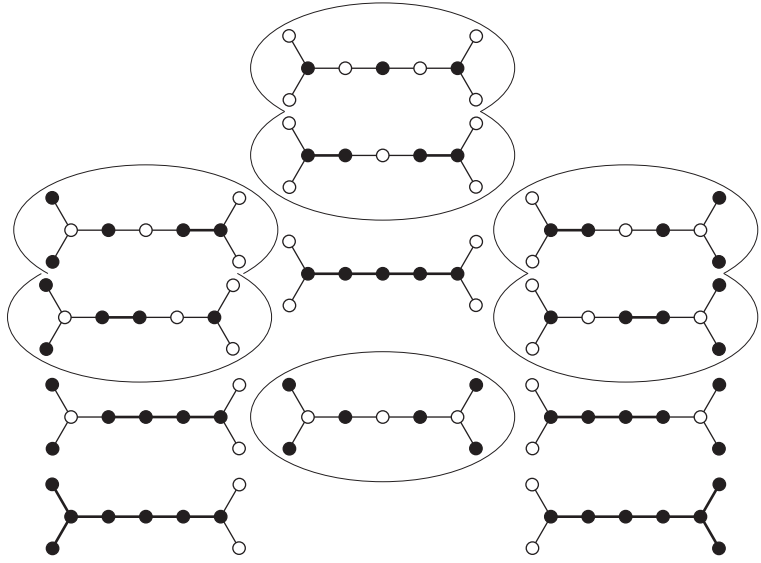}
	\caption{Projective layering of the Veldkamp space of  $\mathcal{C}(\widetilde{D}_8)$.}
	\label{extD8-fig2}
\end{figure}

\begin{figure}[pth!]
	\centering
	\centerline{\includegraphics[width=7.0truecm]{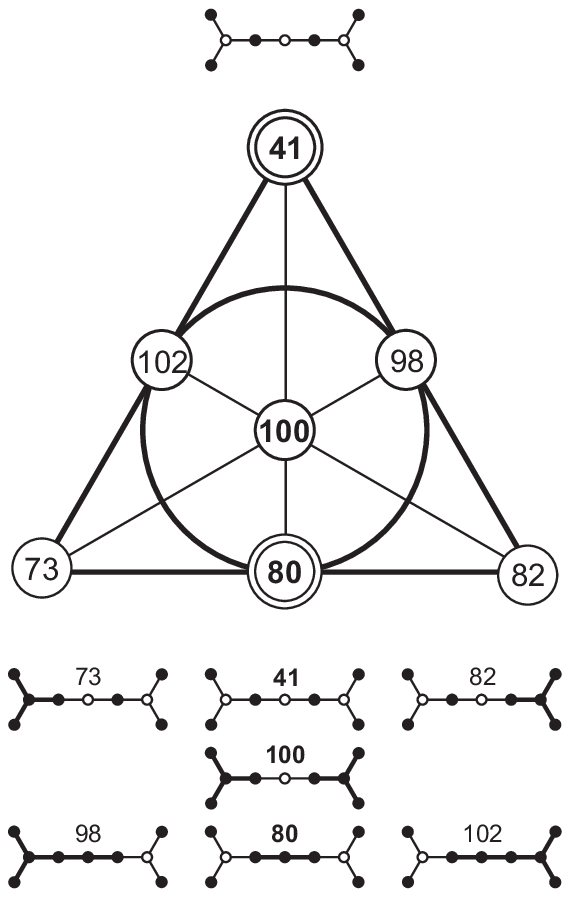} \includegraphics[width=7.0truecm]{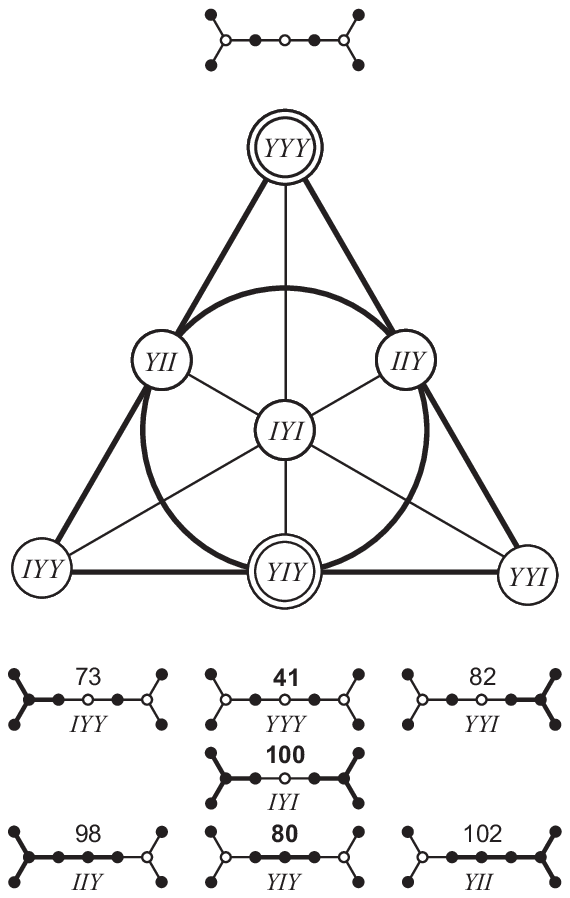}}
	\caption{The explicit structure of the unique Fano plane whose two points are represented by the two exceptional hyperplanes of $\mathcal{C}(\widetilde{D}_8)$ and the associated three-qubit labelling of its points. The four heavy lines and the six points on them form a Pasch configuration.}
	\label{extD8-fig5}
\end{figure}

This is the last case to be dealt with in some detail. Using a computer, we have found that the Veldkamp space of $\mathcal{C}(\widetilde{D}_8)$ has 105 points and 876 size-three lines, exhibiting projective layering as shown in Figure \ref{extD8-fig2}.
From the figure we find out that this Veldkamp space includes one PG$(5,2)$ (1st row) and one PG$(4,2)$ (2nd row), the two having the distinguished PG$(3,2)$ in common (3rd row, middle).
Next, we have here other four PG$(3,2)$s (3rd row, left- and right-hand side), forming two complementary pairs. 
As before, there are two special, disjoint lines in the distinguished PG$(3,2)$, which the latter shares with either of PG$(3,2)$s in both complementary pairs. However, the most interesting object for us is here the Fano plane represented by the hyperplane depicted in the middle of the 4th row of Figure \ref{extD8-fig2}. This hyperplane, as well as the one whose only additional point is point 4, are two exceptional hyperplanes. The corresponding Fano plane, together with all the seven hyperplanes\footnote{The numbering of hyperplanes follows here the same scheme/procedure as adopted in the $n=7$ case.}  that represent its points, are illustrated
 in Figure \ref{extD8-fig5}, {\it left}. Here, the two exceptional hyperplanes are $H_{41}$ and $H_{80}$. They lie on the common Veldkamp line whose third point is $H_{100}$.  Disregarding this point and all the three Veldkamp lines passing through it, we are left with a {\it natural} copy of the so-called {\it Pasch} configuration (thick lines) -- the unique point-line incidence geometry of six points and four lines, with two lines through a point and three points on a line \cite{pasch}.  

\begin{figure}[t]
	\centering
	\centerline{\includegraphics[width=6.0truecm]{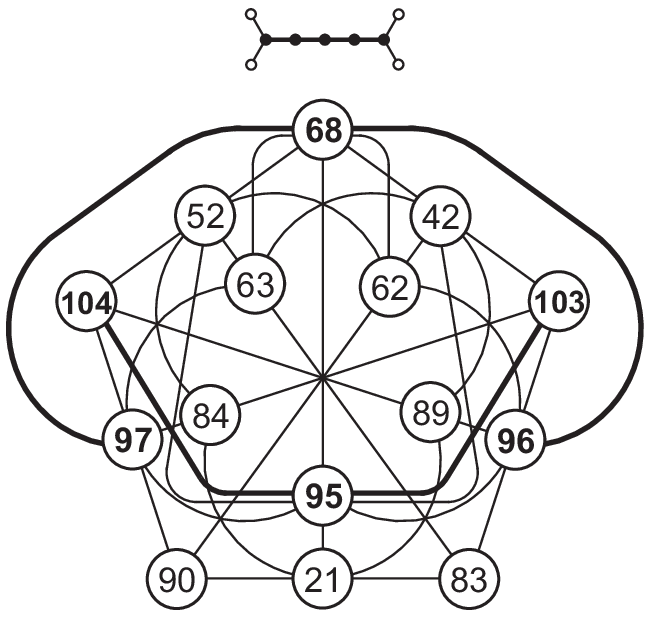}\hspace*{0.5cm}\includegraphics[width=6.0truecm]{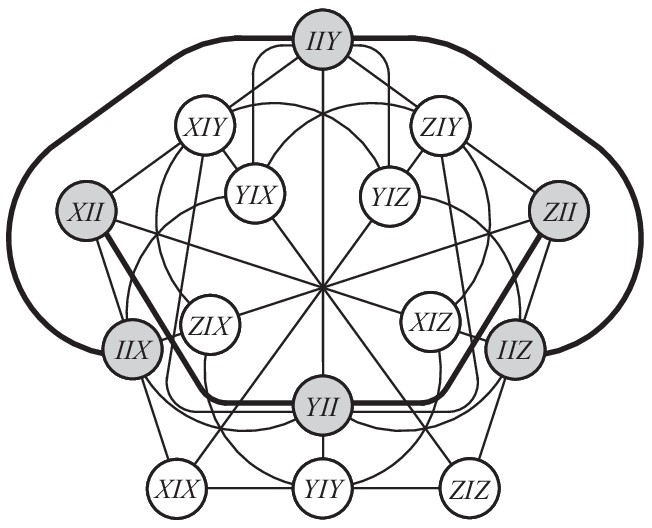}}
	\caption{The hyperplane composition ({\it left}) and its three-qubit counterpart ({\it right}) of the distinguished PG$(3,2)$; points of the two special lines are highlighted in boldface and their corresponding group elements are shaded.}
	\label{extD8-fig3bc}
\end{figure}

This case is remarkable in that the most natural labeling of the (nine) vertices of $\widetilde{D}_8$ employs elements of the {\it three}-qubit Pauli group, in particular
$$0 \rightarrow XII,~ 1 \rightarrow ZII,~ 2 \rightarrow YII,~ 3 \rightarrow IXI,~ 4 \rightarrow IYI,~ 5 \rightarrow IZI,~ 6 \rightarrow IIY,~7 \rightarrow IIX,~8 \rightarrow IIZ.$$ It represents no difficulty to verify that this labelling yields a one-to-one correspondence between 63 elements of the three-qubit Pauli group and 63 points of the PG$(5,2)$. 
Under this correspondence, our distinguished PG$(3,2)$, whose composition is depicted in 
Figure \ref{extD8-fig3bc}, {\it left}, acquires the three-qubit lettering as shown in Figure \ref{extD8-fig3bc}, {\it right}.  
One explicitly sees a bijection between 15 points of this PG$(3,2)$ and  15 elements of a two-qubit subgroup of the three-qubit Pauli group, the geometry of the subgroup encoded in the selected copy of $W(3,2)$ and a three-qubit version of the Mermin-Peres magic square.
An intringuing feature is also absence of the Pauli matrices $X$ and $Z$ in the three-qubit labels of the points of the `exceptional' Fano plane,  as demonstrated by Figure \ref{extD8-fig5}, {\it right}.

\section{Conclusion}
The aim of the paper was to provide both a finite geometer and a mathematical physicist with new examples of the relevance of the concept of Veldkamp space for their respective fields of research.
In the former case, to show how some well-known diagrams, namely extended Dynkin diagrams $\widetilde{D}_n$ ($4 \leq n$), are through Veldkamp spaces intricately related to interesting hierarchies of binary projective spaces.
In the latter case, to illustrate how the structure of the two-qubit (or three-qubit) Pauli group can naturally be invoked to fill up
such `Veldkamp' relations with interesting physical contents.

There are at least two promising ways for further explorations along the lines outlined in the preceding section. One way is to keep going to higher $n$ and focus, for example, on gradually increasing complexity of the hierarchy of projective spaces generated by exceptional geometric hyperplanes and how this complexity can be expressed in the language of multi-qubit Pauli groups.  In this respect the next interesting cases are $n=10$ and $n=11$. In the former case, the `exceptional' PG$(3,2)$ occurs for the first time and it will be interesting to see how this space differs from the always-present distinguished one in terms of three-qubit labeling(s).
In the latter case, the most natural labeling of the vertices of the Dynkin diagram seems to be furnished by elements of the four-qubit Pauli group, whose symplectic polar space $W(7,2)$ features a number of physically relevant finite-geometrical structures (see, e.\,g., \cite{4qov}).  The other way is to go back to the already addressed cases and, following the pioneering work of one of the authors \cite{gq42}, look at properties of those parts of Veldkamp spaces that are composed of lines having {\it more} than three points.

\section*{Acknowledgments}
This work was supported by the Slovak VEGA Grant Agency, Project 2/0003/16, as well as by the French Conseil R\'egional Research Project RECH-MOB15-000007. We are extremely grateful to  J\'er\^ome Boulmier and Benoit Courtier for their computer assistance.

\end{document}